\documentclass[11pt,twoside]{amsart}

\usepackage{latexsym}
\usepackage{amssymb}

\long\def\ig#1{\relax}
\ig{Thanks to Roberto Minio for this def'n.  Compare the def'n of
\comment in AMSTeX.}
 
\newcount \coefa
\newcount \coefb
\newcount \coefc
\newdimen\tempdimen
\newdimen\xlen
\newdimen\ylen
 
\newcount\tempcounta
\newcount\tempcountb
\newcount\tempcountc
\newcount\tempcountd
\newcount\tempcounte
\newcount\tempcountf
\newcount\xext
\newcount\yext
\newcount\xoff
\newcount\yoff
\newcount\gap%
\newcount\arrowtypea
\newcount\arrowtypeb
\newcount\arrowtypec
\newcount\arrowtyped
\newcount\arrowtypee
\newcount\height
\newcount\width
\newcount\xpos
\newcount\ypos
\newcount\run
\newcount\rise
\newcount\arrowlength
\newcount\halflength
\newcount\arrowtype
\newsavebox{\tempboxa}%
\newsavebox{\tempboxb}%
\newsavebox{\tempboxc}%

\catcode`@=11 
\def\settoheight#1#2{\setbox\@tempboxa\hbox{#2}#1\ht\@tempboxa\relax}%
\def\settodepth#1#2{\setbox\@tempboxa\hbox{#2}#1\dp\@tempboxa\relax}%
\let\ifnextchar=\@ifnextchar
\catcode`@=12 
 
\def\putbox(#1,#2)#3{\put(#1,#2){\makebox(0,0){#3}}}

\def\setsqparms[#1`#2`#3`#4;#5`#6]{%
\settripairparms[#1`#2`#3`#4`1;#6]%
\width #5
}
 
\def\settriparms[#1`#2`#3;#4]{\settripairparms[#1`#2`#3`1`1;#4]}%

\def\settripairparms[#1`#2`#3`#4`#5;#6]{%
\arrowtypea #1
\arrowtypeb #2
\arrowtypec #3
\arrowtyped #4
\arrowtypee #5
\height #6
\width #6
}
 
\def\mvector(#1,#2)#3{
\put(0,0){\vector(#1,#2){#3}}%
\put(0,0){\vector(#1,#2){30}}%
}
\def\evector(#1,#2)#3{{
\arrowlength #3
\put(0,0){\vector(#1,#2){\arrowlength}}%
\advance \arrowlength by-30
\put(0,0){\vector(#1,#2){\arrowlength}}%
}}

\def\horsize#1#2{%
\settowidth{\tempdimen}{$#2$}%
#1=\tempdimen
\divide #1 by\unitlength
}
 
\def\vertsize#1#2{%
\settoheight{\tempdimen}{$#2$}%
#1=\tempdimen
\settodepth{\tempdimen}{$#2$}%
\advance #1 by\tempdimen
\divide #1 by\unitlength
}

\def\vertadjust[#1`#2`#3]{%
\vertsize{\tempcounta}{#1}%
\vertsize{\tempcountb}{#2}%
\ifnum \tempcounta<\tempcountb \tempcounta=\tempcountb \fi
\divide\tempcounta by2
\vertsize{\tempcountb}{#3}%
\ifnum \tempcountb>0 \advance \tempcountb by20 \fi
\ifnum \tempcounta<\tempcountb \tempcounta=\tempcountb \fi
}
 
\def\horadjust[#1`#2`#3]{%
\horsize{\tempcounta}{#1}%
\horsize{\tempcountb}{#2}%
\ifnum \tempcounta<\tempcountb \tempcounta=\tempcountb \fi
\divide\tempcounta by20
\horsize{\tempcountb}{#3}%
\ifnum \tempcountb>0 \advance \tempcountb by60 \fi
\ifnum \tempcounta<\tempcountb \tempcounta=\tempcountb \fi
}
 
\ig{ In this procedure, #1 is the paramater that sticks out all the way,
#2 sticks out the least and #3 is a label sticking out half way.  #4 is
the amount of the offset.}
 
\def\sladjust[#1`#2`#3]#4{%
\tempcountc=#4
\horsize{\tempcounta}{#1}%
\divide \tempcounta by2
\horsize{\tempcountb}{#2}%
\divide \tempcountb by2
\advance \tempcountb by-\tempcountc
\ifnum \tempcounta<\tempcountb \tempcounta=\tempcountb\fi
\divide \tempcountc by2
\horsize{\tempcountb}{#3}%
\advance \tempcountb by-\tempcountc
\ifnum \tempcountb>0 \advance \tempcountb by80\fi
\ifnum \tempcounta<\tempcountb \tempcounta=\tempcountb\fi
\advance\tempcounta by20
}
 
\def\putvector(#1,#2)(#3,#4)#5#6{{%
\xpos=#1
\ypos=#2
\run=#3
\rise=#4
\arrowlength=#5
\arrowtype=#6
\ifnum \arrowtype<0
    \ifnum \run=0
        \advance \ypos by-\arrowlength
    \else
        \tempcounta \arrowlength
        \multiply \tempcounta by\rise
        \divide \tempcounta by\run
        \ifnum\run>0
            \advance \xpos by\arrowlength
            \advance \ypos by\tempcounta
        \else
            \advance \xpos by-\arrowlength
            \advance \ypos by-\tempcounta
        \fi
    \fi
    \multiply \arrowtype by-1
    \multiply \rise by-1
    \multiply \run by-1
\fi
\ifnum \arrowtype=1
    \put(\xpos,\ypos){\vector(\run,\rise){\arrowlength}}%
\else\ifnum \arrowtype=2
    \put(\xpos,\ypos){\mvector(\run,\rise)\arrowlength}%
\else\ifnum\arrowtype=3
    \put(\xpos,\ypos){\evector(\run,\rise){\arrowlength}}%
\fi\fi\fi
}}
 
\def\bfig{\begin{picture}(\xext,\yext)(\xoff,\yoff)}
\def\efig{\end{picture}}

\def\putsplitvector(#1,#2)#3#4{
\xpos #1
\ypos #2
\arrowtype #4
\halflength #3
\arrowlength #3
\gap 140
\advance \halflength by-\gap
\divide \halflength by2
\ifnum \arrowtype=1
    \put(\xpos,\ypos){\line(0,-1){\halflength}}%
    \advance\ypos by-\halflength
    \advance\ypos by-\gap
    \put(\xpos,\ypos){\vector(0,-1){\halflength}}%
\else\ifnum \arrowtype=2
    \put(\xpos,\ypos){\line(0,-1)\halflength}%
    \put(\xpos,\ypos){\vector(0,-1)3}%
    \advance\ypos by-\halflength
    \advance\ypos by-\gap
    \put(\xpos,\ypos){\vector(0,-1){\halflength}}%
\else\ifnum\arrowtype=3
    \put(\xpos,\ypos){\line(0,-1)\halflength}%
    \advance\ypos by-\halflength
    \advance\ypos by-\gap
    \put(\xpos,\ypos){\evector(0,-1){\halflength}}%
\else\ifnum \arrowtype=-1
    \advance \ypos by-\arrowlength
    \put(\xpos,\ypos){\line(0,1){\halflength}}%
    \advance\ypos by\halflength
    \advance\ypos by\gap
    \put(\xpos,\ypos){\vector(0,1){\halflength}}%
\else\ifnum \arrowtype=-2
    \advance \ypos by-\arrowlength
    \put(\xpos,\ypos){\line(0,1)\halflength}%
    \put(\xpos,\ypos){\vector(0,1)3}%
    \advance\ypos by\halflength
    \advance\ypos by\gap
    \put(\xpos,\ypos){\vector(0,1){\halflength}}%
\else\ifnum\arrowtype=-3
    \advance \ypos by-\arrowlength
    \put(\xpos,\ypos){\line(0,1)\halflength}%
    \advance\ypos by\halflength
    \advance\ypos by\gap
    \put(\xpos,\ypos){\evector(0,1){\halflength}}%
\fi\fi\fi\fi\fi\fi
}
 
\def\setpos(#1,#2){\xpos=#1 \ypos#2}
 
\def\putmorphism(#1)(#2,#3)[#4`#5`#6]#7#8#9{{%
\run #2
\rise #3
\ifnum\rise=0
  \puthmorphism(#1)[#4`#5`#6]{#7}{#8}{#9}%
\else\ifnum\run=0
  \putvmorphism(#1)[#4`#5`#6]{#7}{#8}{#9}%
\else
\setpos(#1)%
\arrowlength #7
\arrowtype #8
\ifnum\run=0
\else\ifnum\rise=0
\else
\ifnum\run>0
    \coefa=1
\else
   \coefa=-1
\fi
\ifnum\arrowtype>0
   \coefb=0
   \coefc=-1
\else
   \coefb=\coefa
   \coefc=1
   \arrowtype=-\arrowtype
\fi
\width=2
\multiply \width by\run
\divide \width by\rise
\ifnum \width<0  \width=-\width\fi
\advance\width by60
\if l#9 \width=-\width\fi
\putbox(\xpos,\ypos){$#4$}
{\multiply \coefa by\arrowlength
\advance\xpos by\coefa
\multiply \coefa by\rise
\divide \coefa by\run
\advance \ypos by\coefa
\putbox(\xpos,\ypos){$#5$} }%
{\multiply \coefa by\arrowlength
\divide \coefa by2
\advance \xpos by\coefa
\advance \xpos by\width
\multiply \coefa by\rise
\divide \coefa by\run
\advance \ypos by\coefa
\if l#9%
   \put(\xpos,\ypos){\makebox(0,0)[r]{$#6$}}%
\else\if r#9%
   \put(\xpos,\ypos){\makebox(0,0)[l]{$#6$}}%
\fi\fi }%
{\multiply \rise by-\coefc
\multiply \run by-\coefc
\multiply \coefb by\arrowlength
\advance \xpos by\coefb
\multiply \coefb by\rise
\divide \coefb by\run
\advance \ypos by\coefb
\multiply \coefc by70
\advance \ypos by\coefc
\multiply \coefc by\run
\divide \coefc by\rise
\advance \xpos by\coefc
\multiply \coefa by140
\multiply \coefa by\run
\divide \coefa by\rise
\advance \arrowlength by\coefa
\ifnum \arrowtype=1
   \put(\xpos,\ypos){\vector(\run,\rise){\arrowlength}}%
\else\ifnum\arrowtype=2
   \put(\xpos,\ypos){\mvector(\run,\rise){\arrowlength}}%
\else\ifnum\arrowtype=3
   \put(\xpos,\ypos){\evector(\run,\rise){\arrowlength}}%
\fi\fi\fi}%
\fi\fi
\fi\fi}}

\def\puthmorphism(#1,#2)[#3`#4`#5]#6#7#8{{%
\xpos #1
\ypos #2
\width #6
\arrowlength #6
\putbox(\xpos,\ypos){$#3$\vphantom{$#4$}}%
{\advance \xpos by\arrowlength
\putbox(\xpos,\ypos){\vphantom{$#3$}$#4$}}%
\horsize{\tempcounta}{#3}%
\horsize{\tempcountb}{#4}%
\divide \tempcounta by2
\divide \tempcountb by2
\advance \tempcounta by30
\advance \tempcountb by30
\advance \xpos by\tempcounta
\advance \arrowlength by-\tempcounta
\advance \arrowlength by-\tempcountb
\putvector(\xpos,\ypos)(1,0){\arrowlength}{#7}%
\divide \arrowlength by2
\advance \xpos by\arrowlength
\vertsize{\tempcounta}{#5}%
\divide\tempcounta by2
\advance \tempcounta by20
\if a#8 %
   \advance \ypos by\tempcounta
   \put(\xpos,\ypos){\makebox(0,0){$#5$}}%
\else
   \advance \ypos by-\tempcounta
   \put(\xpos,\ypos){\makebox(0,0){$#5$}}%
\fi
}}
 
\def\putvmorphism(#1,#2)[#3`#4`#5]#6#7#8{{%
\xpos #1
\ypos #2
\arrowlength #6
\arrowtype #7
\settowidth{\xlen}{$#5$}%
\putbox(\xpos,\ypos){$#3$}%
{\advance \ypos by-\arrowlength
\putbox(\xpos,\ypos){$#4$}}%
{\advance\arrowlength by-140
\advance \ypos by-70
\ifdim\xlen>0pt
   \if m#8%
      \putsplitvector(\xpos,\ypos){\arrowlength}{\arrowtype}%
   \else
      \putvector(\xpos,\ypos)(0,-1){\arrowlength}{\arrowtype}%
   \fi
\else
   \putvector(\xpos,\ypos)(0,-1){\arrowlength}{\arrowtype}%
\fi}%
\ifdim\xlen>0pt
   \divide \arrowlength by2
   \advance\ypos by-\arrowlength
   \if l#8%
      \advance \xpos by-40
      \put(\xpos,\ypos){\makebox(0,0)[r]{$#5$}}%
   \else\if r#8%
      \advance \xpos by40
      \put(\xpos,\ypos){\makebox(0,0)[l]{$#5$}}%
   \else
      \putbox(\xpos,\ypos){$#5$}%
   \fi\fi
\fi
}}
 
\def\topadjust[#1`#2`#3]{%
\yoff=10
\vertadjust[#1`#2`{#3}]%
\advance \yext by\tempcounta
\advance \yext by 10
}
 
\def\botadjust[#1`#2`#3]{%
\vertadjust[#1`#2`{#3}]%
\advance \yext by\tempcounta
\advance \yoff by-\tempcounta
}
 
\def\leftadjust[#1`#2`#3]{%
\xoff=0
\horadjust[#1`#2`{#3}]%
\advance \xext by\tempcounta
\advance \xoff by-\tempcounta
}
 
\def\rightadjust[#1`#2`#3]{%
\horadjust[#1`#2`{#3}]%
\advance \xext by\tempcounta
}
 
\def\rightsladjust[#1`#2`#3]{%
\sladjust[#1`#2`{#3}]{\width}%
\advance \xext by\tempcounta
}
 
\def\leftsladjust[#1`#2`#3]{%
\xoff=0
\sladjust[#1`#2`{#3}]{\width}%
\advance \xext by\tempcounta
\advance \xoff by-\tempcounta
}
 
\def\adjust[#1`#2;#3`#4;#5`#6;#7`#8]{%
\topadjust[#1``{#2}]
\leftadjust[#3``{#4}]
\rightadjust[#5``{#6}]
\botadjust[#7``{#8}]}

\def\putsquare(#1)[#2`#3`#4`#5;#6`#7`#8`#9]{%
\setpos(#1)
\puthmorphism(\xpos,\ypos)[#4`#5`{#9}]{\width}{\arrowtyped}b%
\advance\ypos by \height
\puthmorphism(\xpos,\ypos)[#2`#3`{#6}]{\width}{\arrowtypea}a%
\putvmorphism(\xpos,\ypos)[``{#7}]{\height}{\arrowtypeb}l%
\advance\xpos by \width
\putvmorphism(\xpos,\ypos)[``{#8}]{\height}{\arrowtypec}r%
}
 
\def\square[#1`#2`#3`#4;#5`#6`#7`#8]{{
\xext=\width                              
\yext=\height                             
\topadjust[#1`#2`{#5}]
\botadjust[#3`#4`{#8}]
\leftadjust[#1`#3`{#6}]
\rightadjust[#2`#4`{#7}]
\begin{picture}(\xext,\yext)(\xoff,\yoff)
\putsquare(0,0)[#1`#2`#3`#4;#5`#6`#7`{#8}]
\end{picture}
}}

\def\putptriangle(#1,#2)[#3`#4`#5;#6`#7`#8]{%
\xpos=#1 \ypos=#2
\advance\ypos by \height
\puthmorphism(\xpos,\ypos)[#3`#4`{#6}]{\height}{\arrowtypea}a%
\putvmorphism(\xpos,\ypos)[`#5`{#7}]{\height}{\arrowtypeb}l%
\advance\xpos by\height
\putmorphism(\xpos,\ypos)(-1,-1)[``{#8}]{\height}{\arrowtypec}r%
}
 
\def\ptriangle[#1`#2`#3;#4`#5`#6]{{
\width=\height                         
\xext=\width                           
\yext=\width                           
\topadjust[#1`#2`{#4}]
\botadjust[#3``]
\leftadjust[#1`#3`{#5}]
\rightsladjust[#2`#3`{#6}]
\begin{picture}(\xext,\yext)(\xoff,\yoff)
\putptriangle(0,0)[#1`#2`#3;#4`#5`{#6}]%
\end{picture}%
}}

\def\putqtriangle(#1,#2)[#3`#4`#5;#6`#7`#8]{%
\xpos=#1 \ypos=#2
\advance\ypos by\height
\puthmorphism(\xpos,\ypos)[#3`#4`{#6}]{\height}{\arrowtypea}a%
\putmorphism(\xpos,\ypos)(1,-1)[``{#7}]{\height}{\arrowtypeb}l%
\advance\xpos by\height
\putvmorphism(\xpos,\ypos)[`#5`{#8}]{\height}{\arrowtypec}r%
}
 
\def\qtriangle[#1`#2`#3;#4`#5`#6]{{
\width=\height                         
\xext=\width                           
\yext=\height                          
\topadjust[#1`#2`{#4}]
\botadjust[#3``]
\leftsladjust[#1`#3`{#5}]
\rightadjust[#2`#3`{#6}]
\begin{picture}(\xext,\yext)(\xoff,\yoff)
\putqtriangle(0,0)[#1`#2`#3;#4`#5`{#6}]%
\end{picture}%
}}

\def\putdtriangle(#1,#2)[#3`#4`#5;#6`#7`#8]{%
\xpos=#1 \ypos=#2
\puthmorphism(\xpos,\ypos)[#4`#5`{#8}]{\height}{\arrowtypec}b%
\advance\xpos by \height \advance\ypos by\height
\putmorphism(\xpos,\ypos)(-1,-1)[``{#6}]{\height}{\arrowtypea}l%
\putvmorphism(\xpos,\ypos)[#3``{#7}]{\height}{\arrowtypeb}r%
}
 
\def\dtriangle[#1`#2`#3;#4`#5`#6]{{
\width=\height                         
\xext=\width                           
\yext=\height                          
\topadjust[#1``]
\botadjust[#2`#3`{#6}]
\leftsladjust[#2`#1`{#4}]
\rightadjust[#1`#3`{#5}]
\begin{picture}(\xext,\yext)(\xoff,\yoff)
\putdtriangle(0,0)[#1`#2`#3;#4`#5`{#6}]%
\end{picture}%
}}

\def\putbtriangle(#1,#2)[#3`#4`#5;#6`#7`#8]{%
\xpos=#1 \ypos=#2
\puthmorphism(\xpos,\ypos)[#4`#5`{#8}]{\height}{\arrowtypec}b%
\advance\ypos by\height
\putmorphism(\xpos,\ypos)(1,-1)[``{#7}]{\height}{\arrowtypeb}r%
\putvmorphism(\xpos,\ypos)[#3``{#6}]{\height}{\arrowtypea}l%
}
 
\def\btriangle[#1`#2`#3;#4`#5`#6]{{
\width=\height                         
\xext=\width                           
\yext=\height                          
\topadjust[#1``]
\botadjust[#2`#3`{#6}]
\leftadjust[#1`#2`{#4}]
\rightsladjust[#3`#1`{#5}]
\begin{picture}(\xext,\yext)(\xoff,\yoff)
\putbtriangle(0,0)[#1`#2`#3;#4`#5`{#6}]%
\end{picture}%
}}

\def\putAtriangle(#1,#2)[#3`#4`#5;#6`#7`#8]{%
\xpos=#1 \ypos=#2
{\multiply \height by2
\puthmorphism(\xpos,\ypos)[#4`#5`{#8}]{\height}{\arrowtypec}b}%
\advance\xpos by\height \advance\ypos by\height
\putmorphism(\xpos,\ypos)(-1,-1)[#3``{#6}]{\height}{\arrowtypea}l%
\putmorphism(\xpos,\ypos)(1,-1)[``{#7}]{\height}{\arrowtypeb}r%
}
 
\def\Atriangle[#1`#2`#3;#4`#5`#6]{{
\width=\height                         
\xext=\width                           
\yext=\height                          
\topadjust[#1``]
\botadjust[#2`#3`{#6}]
\multiply \xext by2 
\leftsladjust[#2`#1`{#4}]
\rightsladjust[#3`#1`{#5}]
\begin{picture}(\xext,\yext)(\xoff,\yoff)%
\putAtriangle(0,0)[#1`#2`#3;#4`#5`{#6}]%
\end{picture}%
}}

\def\putAtrianglepair(#1,#2)[#3]{\xpos=#1 \ypos=#2%
\putAtrianglepairx[#3]}
\def\putAtrianglepairx[#1`#2`#3`#4;#5`#6`#7`#8`#9]{%
\puthmorphism(\xpos,\ypos)[#2`#3`{#8}]{\height}{\arrowtyped}b%
\advance\xpos by\height
\puthmorphism(\xpos,\ypos)[\phantom{#3}`#4`{#9}]{\height}{\arrowtypee}b%
\advance\ypos by\height
\putmorphism(\xpos,\ypos)(-1,-1)[#1``{#5}]{\height}{\arrowtypea}l%
\putvmorphism(\xpos,\ypos)[``{#6}]{\height}{\arrowtypeb}m%
\putmorphism(\xpos,\ypos)(1,-1)[``{#7}]{\height}{\arrowtypec}r%
}
 
\def\Atrianglepair[#1`#2`#3`#4;#5`#6`#7`#8`#9]{{%
\width=\height
\xext=\width
\yext=\height
\topadjust[#1``]%
\vertadjust[#2`#3`{#8}]
\tempcountd=\tempcounta                       
\vertadjust[#3`#4`{#9}]
\ifnum\tempcounta<\tempcountd                 
\tempcounta=\tempcountd\fi                    
\advance \yext by\tempcounta                  
\advance \yoff by-\tempcounta                 
\multiply \xext by2 
\leftsladjust[#2`#1`{#5}]
\rightsladjust[#4`#1`{#7}]%
\begin{picture}(\xext,\yext)(\xoff,\yoff)%
\putAtrianglepair(0,0)[#1`#2`#3`#4;#5`#6`#7`#8`{#9}]%
\end{picture}%
}}

\def\putVtriangle(#1,#2)[#3`#4`#5;#6`#7`#8]{%
\xpos=#1 \ypos=#2
\advance\ypos by\height
{\multiply\height by2
\puthmorphism(\xpos,\ypos)[#3`#4`{#6}]{\height}{\arrowtypea}a}%
\putmorphism(\xpos,\ypos)(1,-1)[`#5`{#7}]{\height}{\arrowtypeb}l%
\advance\xpos by\height
\advance\xpos by\height
\putmorphism(\xpos,\ypos)(-1,-1)[``{#8}]{\height}{\arrowtypec}r%
}
 
\def\Vtriangle[#1`#2`#3;#4`#5`#6]{{
\width=\height                         
\xext=\width                           
\yext=\height                          
\topadjust[#1`#2`{#4}]
\botadjust[#3``]
\multiply \xext by2 
\leftsladjust[#1`#3`{#5}]
\rightsladjust[#2`#3`{#6}]
\begin{picture}(\xext,\yext)(\xoff,\yoff)%
\putVtriangle(0,0)[#1`#2`#3;#4`#5`{#6}]%
\end{picture}%
}}

\def\putVtrianglepair(#1,#2)[#3]{\xpos=#1 \ypos=#2%
\putVtrianglepairx[#3]}
\def\putVtrianglepairx[#1`#2`#3`#4;#5`#6`#7`#8`#9]{%
\advance\ypos by\height
\putmorphism(\xpos,\ypos)(1,-1)[`#4`{#7}]{\height}{\arrowtypec}l%
\puthmorphism(\xpos,\ypos)[#1`#2`{#5}]{\height}{\arrowtypea}a%
\advance\xpos by\height
\puthmorphism(\xpos,\ypos)[\phantom{#2}`#3`{#6}]{\height}{\arrowtypeb}a%
\putvmorphism(\xpos,\ypos)[``{#8}]{\height}{\arrowtyped}m%
\advance\xpos by\height
\putmorphism(\xpos,\ypos)(-1,-1)[``{#9}]{\height}{\arrowtypee}r%
}
 
\def\Vtrianglepair[#1`#2`#3`#4;#5`#6`#7`#8`#9]{{%
\xoff=0
\yoff=2 
\xext=\height                  
\width=\height                 
\yext=\height                  
\vertadjust[#1`#2`{#5}]
\tempcountd=\tempcounta        
\vertadjust[#2`#3`{#6}]
\ifnum\tempcounta<\tempcountd
\tempcounta=\tempcountd\fi
\advance \yext by\tempcounta
\botadjust[#4``]%
\multiply \xext by2
\leftsladjust[#1`#4`{#7}]%
\rightsladjust[#3`#4`{#9}]%
\begin{picture}(\xext,\yext)(\xoff,\yoff)%
\putVtrianglepair(0,0)[#1`#2`#3`#4;#5`#6`#7`#8`{#9}]%
\end{picture}%
}}

\def\putCtriangle(#1,#2)[#3`#4`#5;#6`#7`#8]{%
\xpos=#1 \ypos=#2
\advance\ypos by\height
\putmorphism(\xpos,\ypos)(1,-1)[``{#8}]{\height}{\arrowtypec}l%
\advance\xpos by\height
\advance\ypos by\height
\putmorphism(\xpos,\ypos)(-1,-1)[#3`#4`{#6}]{\height}{\arrowtypea}l%
{\multiply\height by 2
\putvmorphism(\xpos,\ypos)[`#5`{#7}]{\height}{\arrowtypeb}r}%
}
 
\def\Ctriangle[#1`#2`#3;#4`#5`#6]{{
\width=\height                          
\xext=\width                            
\yext=\height                           
\multiply \yext by2 
\topadjust[#1``]
\botadjust[#3``]
\sladjust[#2`#1`{#4}]{\width}
\tempcountd=\tempcounta                 
\sladjust[#2`#3`{#6}]{\width}
\ifnum \tempcounta<\tempcountd          
\tempcounta=\tempcountd\fi              
\advance \xext by\tempcounta            
\advance \xoff by-\tempcounta           
\rightadjust[#1`#3`{#5}]
\begin{picture}(\xext,\yext)(\xoff,\yoff)%
\putCtriangle(0,0)[#1`#2`#3;#4`#5`{#6}]%
\end{picture}%
}}

\def\putDtriangle(#1,#2)[#3`#4`#5;#6`#7`#8]{%
\xpos=#1 \ypos=#2
\advance\xpos by\height \advance\ypos by\height
\putmorphism(\xpos,\ypos)(-1,-1)[``{#8}]{\height}{\arrowtypec}r%
\advance\xpos by-\height \advance\ypos by\height
\putmorphism(\xpos,\ypos)(1,-1)[`#4`{#7}]{\height}{\arrowtypeb}r%
{\multiply\height by 2
\putvmorphism(\xpos,\ypos)[#3`#5`{#6}]{\height}{\arrowtypea}l}%
}
 
\def\Dtriangle[#1`#2`#3;#4`#5`#6]{{
\width=\height                         
\xext=\height                          
\yext=\height                          
\multiply \yext by2 
\topadjust[#1``]
\botadjust[#3``]
\leftadjust[#1`#3`{#4}]
\sladjust[#2`#1`{#4}]{\height}
\tempcountd=\tempcountd                
\sladjust[#2`#3`{#6}]{\height}
\ifnum \tempcounta<\tempcountd         
\tempcounta=\tempcountd\fi             
\advance \xext by\tempcounta           
\begin{picture}(\xext,\yext)(\xoff,\yoff)
\putDtriangle(0,0)[#1`#2`#3;#4`#5`{#6}]%
\end{picture}%
}}

\def\setrecparms[#1`#2]{\width=#1 \height=#2}%
\def\recurse[#1`#2`#3`#4;#5`#6`#7`#8`#9]{{%
\settowidth{\tempdimen}{#1}
\ifdim\tempdimen=0pt
  \savebox{\tempboxa}{\hbox{#2}}%
  \savebox{\tempboxb}{\hbox{#4}}%
  \savebox{\tempboxc}{\hbox{#7}}%
\else
  \savebox{\tempboxa}{\hbox{$\hbox{#1}\times\hbox{#2}$}}%
  \savebox{\tempboxb}{\hbox{$\hbox{#1}\times\hbox{#4}$}}%
  \savebox{\tempboxc}{\hbox{$\hbox{#1}\times\hbox{#7}$}}%
\fi
\tempcounte=\height
\divide\tempcounte by 2
\tempcountf=\tempcounte
\advance\tempcountf by \width
\xext=\tempcountf \yext=\height
\topadjust[#2`\usebox{\tempboxa}`{#5}]%
\botadjust[#4`\usebox{\tempboxb}`{#9}]%
\sladjust[#3`#2`{#6}]{\tempcounte}%
\tempcountd=\tempcounta
\sladjust[#3`#4`{#8}]{\tempcounte}%
\ifnum \tempcounta<\tempcountd
\tempcounta=\tempcountd\fi
\advance \xext by\tempcounta
\advance \xoff by-\tempcounta
\rightadjust[\usebox{\tempboxa}`\usebox{\tempboxb}`\usebox{\tempboxc}]%
\bfig
{\settriparms[-1`1`1;\tempcounte]%
\putCtriangle(0,0)[`#3`;#6`#7`{#8}]}%
\arrowtypea=-1 \arrowtypeb=0 \arrowtypec=1 \arrowtyped=-1
\putsquare(\tempcounte,0)[#2`\usebox{\tempboxa}`#4`\usebox{\tempboxb};%
#5``\usebox{\tempboxc}`#9]%
\efig
}}

\newcommand{\Id}{\mathrm{Id}} 

\newcommand{\dopu}{{:}\allowbreak\ }
\newcommand{\eps}{\varepsilon}
\newcommand{\de}{\delta}
\newcommand{\cal}{\mathcal}

\newcommand{\Om}{(\Omega, \Sigma, \mu)}
\def\DP{Daugavet property}
\newcommand{\Ainfty}{\AAA_\infty}

\newcommand{\loglike}[1]{\mathop{\rm #1}\nolimits}

\newcommand{\supp}{\loglike{supp}}
\newcommand{\Lin}{\loglike{\overline{lin}}}
\newcommand{\unc}{\loglike{unc}}

\newcommand{\Sig}{\Sigma}
\newcommand{\Sip}{\Sigma^+}

\newcommand{\Ome}{\Omega}

\newcommand{\0}{{\emptyset}}
\newcommand{\N}{{\mathbb N}}

\newcommand{\LL}{\mathcal L}
\newcommand{\KK}{\mathcal K}

\newcommand{\AAA}{\mathcal A}

\theoremstyle{plain}
\newtheorem{thm}{Theorem}[section]

\newtheorem{prop}[thm]{Proposition}
\newtheorem{cor}[thm]{Corollary}
\newtheorem{lemma}[thm]{Lemma}

\theoremstyle{definition}
\newtheorem{definition}[thm]{Definition}

\theoremstyle{remark}

\newtheorem{rem}[thm]{Remark}

\numberwithin{equation}{section}

\newcommand{\rest}[2]{#1\raisebox{-0.3ex}{\mbox{$\mid_{#2}$}}}

\newcommand{\begsta}{\begin{statements}}
\newcommand{\begaeq}{\begin{aequivalenz}}
\def\endsta{\end{statements}}
\def\endaeq{\end{aequivalenz}}
\newcommand{\bea}{\begin{eqnarray*}}
\newcommand{\eea}{\end{eqnarray*}}

\newcounter{abc}   
\newcounter{iiiii} 

\newenvironment{aequivalenz}
{\setcounter{iiiii}{0}
\begin{list}%
{{\rm (\roman{iiiii})}}
{\usecounter{iiiii}
\parsep=0pt plus 1pt
\topsep=1pt plus 2pt minus 1pt
\itemsep=1pt plus 2pt minus 1pt
\leftmargin=3\baselineskip
\labelsep=.6\baselineskip
\labelwidth=2.4\baselineskip
\rightmargin 0pt}%
}%
{\end{list}}

\newenvironment{statements}%
{\setcounter{abc}{0}
\begin{list}%
{{\rm (\alph{abc})}}
{\usecounter{abc}
\parsep=0pt plus 1pt
\topsep=1pt plus 2pt minus 1pt
\itemsep=1pt plus 2pt minus 1pt
\leftmargin=3\baselineskip
\labelsep=.6\baselineskip
\labelwidth=2.4\baselineskip
\rightmargin 0pt}%
}%
{\end{list}}

\begin{document}

\title[Unconditionally convergent Series of operators]%
{Unconditionally convergent series of operators and narrow operators on $L_1$}

\author{Vladimir Kadets, Nigel Kalton and Dirk Werner}


\address{Faculty of Mechanics and Mathematics, Kharkov National
University,\linebreak
 pl.~Svobody~4,  61077~Kharkov, Ukraine}
\email{vova1kadets@yahoo.com}

\curraddr{Department of Mathematics, University of Missouri,
Columbia MO 65211}

\address{Department of Mathematics, University of Missouri,
Columbia MO 65211}
\email{nigel@math.missouri.edu}

\address{Department of Mathematics, Freie Universit\"at Berlin,
Arnimallee~2--6, \qquad {}\linebreak D-14\,195~Berlin, Germany}
\email{werner@math.fu-berlin.de}

\thanks{The work of the first-named author
was supported by a fellowship from the \textit{Alexander-von-Humboldt
Stiftung}. The second-named author was supported by NSF grant DMS-9870027.}



\begin{abstract}
We introduce a class of operators on $L_1$ that
is stable under taking
sums of pointwise unconditionally convergent series,  
contains all compact operators  and 
does not contain isomorphic embeddings. It follows that any operator
from $L_1$ into a space with an unconditional basis belongs to this class.
\end{abstract}

\maketitle

\thispagestyle{empty}


\section{Introduction}

A famous  theorem due to A.~Pe{\l}czy\'nski
\cite{Pel-Imposs} states that $L_{1}[0,1]$ cannot 
be embedded in a space with an unconditional basis. A somewhat stronger version
is also true \cite{KadShv}: If an operator $J\dopu  L_{1}[0,1] \to X$ is 
bounded from below, then 
it cannot be represented as a pointwise unconditionally convergent series 
of compact operators.  This last theorem   in fact also holds  for embedding
operators $J\dopu  E \to X$  if $E$ has the Daugavet property; 
see \cite{KadSSW}.

We wish to rephrase the theorem using the following definition.

\begin{definition} \label{unc}
Let $\mathcal{ U}$ be a linear subspace of $\LL (E,X)$, the space of
bounded linear operators from $E$ into $X$. By
$\unc (\mathcal{ U})$ we denote the set of all operators which can be 
represented by  pointwise unconditionally convergent series 
of operators from $\mathcal{ U}$.
\end{definition}

In  terms of this definition the above theorem says that
 an isomorphic embedding operator $J\dopu  L_{1}[0,1] \to X$
does not belong to $\unc (\KK (L_1[0,1],X))$, where $\KK(E,X)$ stands
 for the space of compact operators from $E$ into $X$. 

Clearly, one can iterate the operation ``$\unc$'' and consider the classes
$$
\unc (\unc (\KK (L_1[0,1],X))),  \quad \unc (\unc (\unc (\KK
(L_1[0,1],X)))),
$$ 
etc.
Thus the  question arises whether one can obtain an 
isomorphic embedding operator through such a chain of iterations;
indeed it is not clear at the outset whether possibly 
$\unc (\unc (\KK (E,X))) = \unc (\KK (E,X))$.

A natural approach to generalise  Pe{\l}czy\'nski's theorem in
this direction is to
find a large  class of operators $T\dopu  L_{1}[0,1] \to X$  
which  is stable under taking
sums of pointwise unconditionally convergent series,  
 contains all compact operators  and 
  does not contain isomorphic embeddings.

It was shown by R.~Shvidkoy in his Ph.D.~Thesis \cite{Shv0}
and independently in \cite{KadPop1} 
that in the case 
$X=L_{1}[0,1]$, the  PP-narrow operators on
$L_{1}[0,1]$ form such a class. 
 Here is the definition.

Let $\Om$ be
 a fixed nonatomic probability space and   $L_p=L_p\Om $.  
By $\Sip$ we denote the collection of all
measurable subsets of $\Ome$ having nonzero measure.

\begin{definition} \label{nar2}
Let $A \in \Sip$. 
\begsta
\item
A function $f \in L_p$ is said to be
a \textit{sign  supported on $A$} if $f=\chi_{B_1} - \chi_{B_2} $,
where $B_1$ and $B_2$ form a partition of $A$ into two measurable
subsets of equal measure. 
\item
An operator 
$T \in {\LL}(L_p,X)$ is said to be \textit{PP-narrow} if for every 
set $A \in \Sip$ and every $\eps > 0$ there is a sign $f$ 
supported on $A$ with $\|Tf\| \le \eps$. 
\endsta
\end{definition}

The concept of a PP-narrow
operator was introduced by Plichko and Popov in \cite{PliPop} under the name
\textit{narrow operator}. We use the term ``PP-narrow''  in order to
distinguish such operators 
from a related concept of a narrow operator from \cite{KadSW2}, where,
incidentally, PP-narrow operators were called $L_1$-narrow. 
It should be noted that PP-narrow operators appear implicitly in
Rosenthal's papers on sign embeddings (e.g., \cite{Ros-emb}), where an
operator on $L_1$ is called sign preserving if it is not PP-narrow.

Obviously, no embedding operator is PP-narrow. On the other hand it
is clear that a compact operator $T$ is PP-narrow. Indeed, let $(r_n)$ be
a Rademacher sequence supported on a set $A\in \Sip$; i.e., the $r_n$
are stochastically independent with respect to the probability space
$(A,\rest{\Sigma}{A}, \mu/\mu(A))$ and 
$\mu({\{r_n=1\}}) = \mu({\{r_n=-1\}}) = \mu(A)/2$. Then $r_n\to0$
weakly and hence $Tr_n\to0$ in norm.
The same argument shows that weakly compact operators 
on~$L_1$ are PP-narrow, since $L_1$ has the Dunford-Pettis property. 

The aim of this paper is to find a  class of operators
with the above properties that works 
for general $X$ rather than just for $X=L_1[0,1]$. For this purpose we shall
introduce the class of hereditarily PP-narrow (for short HPP-narrow) 
operators in Section~2.
We show that they form a linear space of operators (which is false for
PP-narrow operators, at least for $p>1$), and in Section~3 we derive a
factorisation scheme for unconditional sums of such operators.
This enables us to give an example of a Banach space $X$ for which
$\unc (\unc (\KK (X,X))) \neq \unc (\KK (X,X))$
(Theorem~\ref{theo3.3}). In Section~4 we specialise to the case $p=1$
and obtain that a pointwise unconditionally convergent series of
HPP-narrow operators on $L_1$ is HPP-narrow (Theorem~\ref{thm4.3}). As
a result, it follows that no embedding operator is in any of the spaces
$\unc(\dots (\unc (\KK (L_1,X)))) $. A further
consequence is that every operator from $L_1$ into a space with an
unconditional basis is HPP-narrow and in particular PP-narrow; this
implies that $L_1$ does not even sign-embed into a space with an
unconditional basis. These last results are
due to Rosenthal in his unpublished paper \cite{Ros-stop} (not only
is this paper unpublished, as a matter of fact 
it has never been written, as Rosenthal
has pointed out to us).

In this paper we deal with real Banach spaces.


\section{Haar-like systems and hereditarily PP-narrow operators}

We start by introducing some notions that will be used throughout the
paper.

Denote 
$$
\AAA_0=\{\0\}, \quad \AAA_n=\{-1,1\}^n, \quad \Ainfty =
\bigcup_{n=0}^\infty \AAA_n. 
$$  
The elements of $\AAA_n$ are $n$-tuples
  of the form $(\alpha_1,\dots,\alpha_n)$ with
$\alpha_k=\pm 1$. For 
$\alpha=(\alpha_1,\dots,\alpha_n) \in \AAA_n$ and 
$\alpha_{n+1} \in \{-1,1\}$ denote by $\alpha,\alpha_{n+1}$ 
the $(n+1)$-tuple
$(\alpha_1,\dots,\alpha_n, \alpha_{n+1})\in \AAA_{n+1}$;
also, put  $\0,\alpha_1 = (\alpha_1)$. The elements of  $\Ainfty$
can be written as a sequence in the following \textit{natural order}:
$$
\0,  \  -1, \  1, \  (-1,-1), \  (-1,1), \ (1,-1), \ (1,1) ,\
(-1,-1,-1), \ \dots .
$$

\begin{definition} \label{haarlike}
Let $A \in \Sip$. 
\begsta
\item
A collection  
$\{A_\alpha \dopu  \alpha \in \Ainfty \} $
of subsets of $A$ is said to be
a \textit{tree of subsets on $A$} 
if $A_\0=A$ and if for every $\alpha \in \Ainfty$ the
subsets $A_{\alpha,1}$ and $A_{\alpha,-1}$
form a partition of $A_\alpha$ into two measurable
subsets of equal measure. 
\item
The collection of functions 
$\{h_\alpha \dopu  \alpha \in \Ainfty \} $
defined by
$h_\alpha = \chi_{A_{\alpha,1} } - \chi_{A_{\alpha,-1}}$
is said to be a \textit{Haar-like system on $A$} (corresponding to
the tree of subsets 
$A_\alpha$, $ \alpha \in \Ainfty$). 
\endsta
\end{definition}

It is easy to see that after deleting the constant function the classical 
Haar system is an example of a Haar-like system. Moreover, every
Haar-like system is equivalent to this example.
In particular we note:

\begin{rem} \label{hl0}
(a) Let 
$\{h_\alpha \dopu  \alpha \in \Ainfty \} $
 be a 
Haar-like system on $A$ corresponding to
a tree of subsets $A_\alpha $, and let 
$1 \le p < \infty$.
Denote by $\Sigma_1$ the $\sigma$-algebra on $A$ generated by 
the subsets $A_\alpha$. Then the system 
$\{h_\alpha \dopu  \alpha \in \Ainfty \} $
 in its natural order
forms a monotone Schauder basis for the subspace $L_p^0(A,\Sigma_1, \mu)$ of
$L_p(A,\Sigma_1, \mu)$ consisting of all $f \in L_p(A,\Sigma_1, \mu)$
with $\int_A f \,d \mu = 0$. 
Note that, for $\alpha \in \AAA_n$, 
$\|h_\alpha\|= \bigl( 2^{-n} \mu(A) \bigr)^{1/p}$ for every
Haar-like system on $A$.

(b) Therefore, if $\eps>0$ and $\{\eps_\alpha\dopu \alpha \in
\Ainfty\} $ is a family of positive numbers such that $\sum_\alpha
\eps_\alpha/\|h_\alpha\| \le \eps/2$ and if $\{x_\alpha\dopu \alpha \in
\Ainfty\} $ is a family of vectors in a Banach space $X$ such that
$\|x_\alpha\|\le \eps_\alpha$, then the mapping $h_\alpha \mapsto
x_\alpha$ extends to a bounded linear operator from
$L_p^0(A,\Sig_1,\mu)$ to $X$ of norm ${\le\eps}$.
\end{rem}

\begin{lemma} \label{hl1}
Let $1 \le p < \infty$ and let $T\dopu L_p \to X$ be a PP-narrow operator. 
\begsta
\item 
For every $A \in \Sip$ and every family of numbers
$\eps_\alpha > 0$ there is a Haar-like system 
$\{h_\alpha \dopu  \alpha \in \Ainfty \} $
on $A$
such that $\|Th_\alpha\| \le \eps_\alpha$ for $\alpha \in \Ainfty $.
\item 
For every 
$\eps > 0$ and every $A \in \Sip$ there is a $\sigma$-algebra 
$\Sigma_\eps \subset \Sigma$ on $A$ such that $(A,\Sigma_\eps, \mu)$
is a nonatomic measure space and the restriction of $T$ to
$L_p^0(A,\Sigma_\eps, \mu)$ has norm  ${\le\eps}$.
\endsta
\end{lemma}

\begin{proof}
To construct a tree of subsets and the corresponding Haar-like
system for (a) we repeatedly apply the
definition of a PP-narrow operator. Namely, let $h_\0$ be a sign
supported on $A$ with $\|Th_\0\| \le \eps_\0$. Put, using the notation
$\{h= x\} = \{\omega\dopu h(\omega)=x\}$, 
$$
A_{-1}= \{ h_\0 =-1\}, \quad  A_{1}= \{h_\0 =1\}. 
$$
Let $h_{-1}$ 
and $h_{1}$ be
signs supported on $A_{-1}$ and $A_{1}$ 
respectively with $\|Th_{\pm 1}\| \allowbreak \le \eps_{\pm1}$; put
\bea
&A_{-1,-1}= \{h_{-1}=-1\},\quad  A_{-1,1}= \{h_{-1}=1\},& \\
&A_{1,-1}= \{h_{1} = -1\},\quad  A_{1,1}=  \{h_{1} =1\}&
\eea
and continue in the above fashion. This yields
 part~(a).

Part~(b)  follows from (a) and
Remark~\ref{hl0}(b).
\end{proof}

For $1 < p < \infty$ the class of PP-narrow operators on $L_p$ is
not stable under taking sums  (see \cite{PliPop}, p.~59);
 this is why we have to
consider a smaller class of operators that we introduce next. 
Incidentally, the stability of PP-narrow operators on $L_1$ under sums
is still an open problem.

\begin{definition} \label{hnar1}
An operator $T\dopu L_p \to X$ is said to be \textit{hereditarily 
PP-narrow}  (\textit{HPP-narrow} for short) if for every $A \in \Sip$ 
and every nonatomic
sub-$\sigma$-algebra $\Sigma_1 \subset \Sigma$ on $A$ the restriction
of $T$ to  $L_p(A,\Sigma_1, \mu)$ is PP-narrow.
\end{definition}

Since every compact operator on $L_p$ is PP-narrow 
and compactness is
inherited by restrictions, compact operators on $L_p$ are HPP-narrow.
On the other hand, the operator
$$
T\dopu L_p([0,1]^2) \to L_p[0,1], \quad 
(Tf)(s) = \int_0^1 f(s,t)\,dt
$$
shows that a PP-narrow operator need not be HPP-narrow.

We now show that the set of HPP-narrow operators forms a subspace of
$\LL(L_p,X)$. 

\begin{prop} \label{prop2.5}
Let $1 \le p < \infty$ and let $U,V \dopu L_p \to X$.  
\begsta
\item \label{hnar1.1}
If $U$ is  PP-narrow and $V$ is  HPP-narrow,
then $U+V$ is  PP-narrow.
\item \label{hnar1.2}
If $U$ and $V$ are both  HPP-narrow,
then $U+V$ is  HPP-narrow as well.
\endsta
\end{prop}

\begin{proof}
(a)  Let  $A \in \Sip$ and $\eps > 0$. By
Lemma~\ref{hl1}(b) there is a $\sigma$-algebra 
$\Sigma_\eps \subset \Sigma$ on $A$ such that $(A,\Sigma_\eps, \mu)$
is a nonatomic measure space and the restriction of $U$ to
$L_p^0(A,\Sigma_\eps, \mu)$ has norm  ${\le\eps}$.
Since $V$ is  HPP-narrow, there is a $\Sigma_\eps$-measurable sign
$f$ supported on $A$ for which $\|Vf\| \le \eps$.
Then $\|(U+V)f\| \le \eps \mu(A)^{1/p} + \eps \le 2\eps$.

(b) follows from (a).
\end{proof}


\section{Unconditionally convergent series of HPP-narrow operators}

In this section we are  going to give an example of a Banach
space $X$ for which 
$$
\Id   \in \unc (\unc (\KK (X,X))) \setminus \unc (\KK(X,X)).
$$

We begin with a factorisation lemma for unconditional 
sums of HPP-narrow operators.

\begin{lemma} \label{lem3.1}
Let $1 \le p < \infty$, $X$ be a Banach space, $T_n \dopu L_p \to X$
be  HPP-narrow operators with $\sum_{n=1}^\infty T_n $ converging
pointwise unconditionally to an operator $T$ and let
$M= \sup_\pm \|\sum_{n=1}^\infty \pm T_n\|$. 
Given $0<\eps <1/2 $,
  there exist a Banach space $Y$ and a factorisation 
\begin{center}
\setlength{\unitlength}{.01em}
\settriparms[1`1`-1;500]
\Vtriangle[L_p`X`Y;T`\tilde{T}`W]
\end{center}
with $\|\tilde{T}\| \le M$, $\|W\| \le 1$, and
 there are a nonatomic sub-$\sigma$-algebra 
$\Sig_1 \subset \Sig$, a Haar-like system $\{h_\alpha \}$
forming a basis for  $L_p^0(\Ome,\Sig_1, \mu)$ and operators 
$U,V \dopu L_p^0(\Ome,\Sig_1, \mu) \to Y$ with 
$U+V=\tilde{T}$ on $L_p^0(\Ome,\Sig_1, \mu)$ such that 
$U$ maps $\{h_\alpha \}$ to a $1$-unconditional basic sequence
and $\|V\| \le \eps$.
\end{lemma}

\begin{proof}
Define $Y$ as the space of all sequences $y=(y_1,y_2, \dots )$, 
$y_n \in X$, such that $\sum_{n=1}^\infty y_n $ converges 
unconditionally in $X$. Equip $Y$ with the natural norm 
$$
\|y\| = \sup_\pm \biggl\|\sum_{n=1}^\infty \pm y_n \biggr\|.
$$
Put $\tilde{T}f = (T_1f, T_2f, \dots )$ and
$Wy= \sum_{n=1}^\infty y_n$. Then $ Y$, $\tilde{T}$ and $W$ satisfy
the desired factorisation scheme.

Our main task is now to define for this  $\tilde{T}$ a Haar-like system 
$\{h_\alpha \}$ and operators $U,V$ as claimed in the lemma. To do this one
 uses a standard blocking technique and the stability of 
HPP-narrow operators under summation (Proposition~\ref{prop2.5}). Namely,
for every $1 \le n < m \le \infty$ define a projection operator 
$P_{n,m} \dopu Y \to Y$ as follows:
$$
P_{n,m}(y_1,y_2, \dots ) = 
(0,0,\dots,0,  y_n, y_{n+1}, \dots, y_{m-1}, 0,0, \dots ).
$$
Let $(\eps_\alpha)$ be positive numbers.
Select an arbitrary sign $h_\0$ supported on $\Ome$  and find 
$n_\0 \in \N$ for which  
$$
\|P_{n_\0, \infty} \tilde{T} h_\0 \| \le \eps_\0.
$$
Put 
$$
U h_\0 = P_{1, n_\0}\tilde{T} h_\0, \quad
Vh_\0 = P_{n_\0, \infty} \tilde{T} h_\0.
$$
The sign $h_\0$ generates a partition of $\Ome$, i.e.,
$$
A_{-1}= \{h_\0 =-1\}, \quad A_{1}= \{h_\0 =1\}.
$$
Since the operator $P_{1, n_\0}\tilde{T}$ is PP-narrow by
Proposition~\ref{prop2.5}, there is 
a sign $h_{-1}$ supported on $A_{-1}$ for which 
$$
\|P_{1, n_\0}\tilde{T}h_{-1}\| \le \frac12 \eps_{-1}.
$$
Find $n_{-1} > n_\0$ such that 
$$
\|P_{n_{-1},\infty}\tilde{T}h_{-1}\| \le \frac12 \eps_{-1}.
$$ 
Put 
$$
U h_{-1} = P_{n_\0, n_{-1}}\tilde{T} h_{-1}, \quad
Vh_{-1} = (P_{1,n_\0}+ P_{n_{-1}, \infty}) \tilde{T} h_{-1}.
$$
Continuing in this fashion we obtain a Haar-like system 
$\{h_\alpha \}$ and operators $U,V \dopu \Lin\{h_\alpha \} \to Y$ 
such that $U+V=\tilde{T}$ on $\Lin\{h_\alpha \}$, $U$ maps 
$\{h_\alpha \}$ to disjoint elements of the sequence space $Y$ 
and hence to a 1-unconditional basic sequence and $V$ maps 
$\{h_\alpha \}$ to elements whose norms are controlled 
by the numbers $\eps_\alpha$; therefore $\|V\|\le \eps$ by
Remark~\ref{hl0}(b) if $\eps_\alpha\to 0$ sufficiently fast. 
\end{proof}

\begin{lemma} \label{lem3.2}
Under the conditions of Lemma~\ref{lem3.1} assume in addition that the 
operator $T$ is bounded from below by a constant $c$; i.e.,
$$
\|Tf\| \ge c \|f\| \qquad\forall f \in L_p.
$$
Then 
$$
M= \sup_\pm \biggl\| \sum_{n=1}^\infty \pm T_n \biggr\| \ge \beta_p c,
$$ 
where $\beta_p$ is the unconditional
constant of the Haar system in $L_p$.
\end{lemma}

\begin{proof}
Let $0<\eps<1/2$.
Under the above conditions the operator $U$ from  Lemma~\ref{lem3.1}
maps a Haar-like system $\{h_\alpha \}$ to a 1-unconditional 
basic sequence. This implies that if $U$ is considered as acting
from $\Lin\{h_\alpha \}$ into $\Lin\{Uh_\alpha \}$, then 
$\|U\| \|U^{-1}\| \ge \beta_p$. On the other hand
$$
\|U\| \le \|\tilde{T}\| + \|V\| \le M + \eps
$$
and 
$$
\|Uf\| \ge \|\tilde{T}f\| -  \eps\|f\|
\ge \|{T}f\| -  \eps\|f\| \ge (c - \eps)\|f\|
$$
for all $f \in \Lin\{h_\alpha \}$, so 
$ \|U^{-1}\| \le (c - \eps)^{-1}$. Hence
we have $(M + \eps)(c - \eps)^{-1}\ge \beta_p$, which yields the desired
inequality since $\eps>0$ was arbitrary. 
\end{proof}

It is known that $\beta_p\to\infty$ if $p\to1$ or $p\to \infty$; in
fact, Burkholder \cite{Burk85} has shown that
$$ 
\beta_p =  \max\Bigl\{ p-1 , \frac1{p-1} \Bigr\} .
$$

\begin{thm} \label{theo3.3}
There exists a Banach space $X$ for which 
$$
\Id \in \unc (\unc (\KK (X,X))) \setminus \unc (\KK(X,X)).
$$
\end{thm}

\begin{proof}
Consider the space $X=L_{p_1} \oplus_2 L_{p_2} \oplus_2 \dots$
where $1<p_n <\infty$ and   $p_n \to 1$.

Suppose that
  $\Id  = \sum_{n=1}^\infty T_n $ 
pointwise unconditionally with compact operators~$T_n $.
The restrictions of $T_n $ to $L_{p_j}$ are also compact
and hence HPP-narrow, so by the previous lemma
$$ 
\sup_\pm \biggl\| \sum_{n=1}^\infty \pm T_n \biggr\| \ge
\sup_\pm \biggl\| \sum_{n=1}^\infty \pm \rest{T_n}{L_{p_j}} \biggr\| 
\ge \beta_{p_j} \to \infty.
$$
So the assumption of pointwise unconditional convergence
of $\sum_{n=1}^\infty T_n$ leads to a contradiction, and hence
$\Id $ does not belong to $\unc (\KK(X,X))$.

On the other hand all the natural projections 
$P_j\dopu  X \to  L_{p_j}$ belong to $\unc (\KK(X,X))$ 
since each $L_{p_j}$ has an unconditional basis.
Taking into account the unconditional representation
$\Id  = \sum_{n=1}^\infty P_n $ we obtain that 
$\Id  \in \unc (\unc (\KK (X,X)))$.
\end{proof}


\section{HPP-narrow operators on $L_{1}$}

In this section we prove the main result of the 
paper, namely that the sum of a pointwise unconditionally
convergent series of HPP-narrow operators on $L_1$ is again an
HPP-narrow operator.

The following lemma implies that the operator $U$ from
Lemma~\ref{lem3.1} factors through~$c_0$.

\begin{lemma} \label{lem4.1}
Let $\{h_\alpha \}$ be  a Haar-like system in $L_1$,
$U\dopu L_1 \to X$ be an operator which maps
$\{h_\alpha \}$ into an unconditional basic sequence.
Then there is a constant $C$ such that for every
element of the form $f = \sum_\alpha a_\alpha h_\alpha$ 
one has
\begin{equation} \label{eq4.1}
\|Uf\| \le C \sup_\alpha|a_\alpha|.
\end{equation}
\end{lemma}

\begin{proof}
Without loss of generality we can assume that 
$\|U\| = 1$, $\|h_\0\| = 1$ and that the unconditional
constant of $\{U h_\alpha \}$ also equals~$1$ 
(one can achieve all these goals by an equivalent renorming of
$X$ and by  multiplication of $\mu$ by a constant).

Let us first remark that for every 
$\alpha=(\alpha_1,\alpha_2, \dots,\alpha_n) \in \AAA_n$
$$
\|\alpha_1 h_\0 + 2\alpha_2h_{\alpha_1} +
4\alpha_3h_{\alpha_1,\alpha_2} + 
\dots + 2^{n-1}\alpha_n h_{\alpha_1, \dots, \alpha_{n-1}}\| \le 2;
$$
indeed, it is easy to check by induction over $n$
that this sum equals 
$$
2^n \chi_{A_{\alpha_1, \dots, \alpha_n}} - \chi_{A_\0}.
$$
Hence
$$
\|\alpha_1 Uh_\0 + 2\alpha_2Uh_{\alpha_1} +
\dots + 2^{n-1}\alpha_n U h_{\alpha_1, \dots, \alpha_{n-1}}\| \le 2,
$$
and, since  $\{U h_\alpha \}$ is a 1-unconditional basic sequence,
$$
\| Uh_\0 + 2Uh_{\alpha_1} +
\dots + 2^{n-1} U h_{\alpha_1, \dots, \alpha_{n-1}}\| \le 2.
$$
Passing  from $n-1$ to $n$ in the last inequality and
averaging over $\alpha \in \AAA_n$ we obtain that
$$
2 \ge 
\biggl\| \frac{1}{2^n} \sum_{\alpha \in \AAA_n}
(Uh_\0 + 2Uh_{\alpha_1} + \dots + 2^{n-1} U h_{\alpha_1, \dots,
  \alpha_n})
\biggr\| 
= \biggl\| \sum _{k=0}^n \sum _{\alpha \in \AAA_k} U h_\alpha \biggr\|.
$$
Again by 1-unconditionality of $\{U h_\alpha \}$  the last 
inequality implies that for all $a_\alpha \in [-1,1]$
$$
 \biggl\|\sum _{k=0}^n 
\sum _{\alpha \in \AAA_k}a_\alpha U h_\alpha \biggr\| \le 2,
$$
which gives (\ref{eq4.1}) with $C=2$.
\end{proof}

An inspection of the proof shows that 
$$
\|Uf\| \le 2\|U\| \beta^2 \sup_\alpha |a_\alpha|
$$
where $\beta$ denotes the unconditional constant of the basic
sequence $(Uh_\alpha)$.

\begin{lemma} \label{lem4.2}
For every Haar-like system  $\{h_\alpha \}$ in $L_1$ supported
on $A$
and every $\de > 0$ there is a sign 
\begin{equation} \label{eq4.2}
f = \sum _{k=0}^\infty
\sum _{\alpha \in \AAA_k}a_\alpha h_\alpha
\end{equation}
supported on $A$ with $\sup_\alpha|a_\alpha| \le \de$.
\end{lemma}

\begin{proof}
Fix an $m \in \N$ such that $1/m \le \de$ and define
$$
f_k = \sum _{\alpha \in \AAA_k}a_\alpha  h_\alpha
$$
as follows: $f_0=\frac1m h_\0$, and for every $\alpha \in \AAA_n$
put $a_\alpha = 1/m$ if 
$| \sum _{k=0}^{n-1}f_k | < 1$ on
$\supp h_\alpha$
and  $a_\alpha = 0$ if 
$| \sum _{k=0}^{n-1}f_k | = 1$ on
$\supp h_\alpha$.
Under this construction all the partial sums of the series
$\sum_{k=0}^\infty f_k$ are bounded by~$1$ in modulus. 
Since $\{f_k\}_{k=0}^\infty$ is
an orthogonal system, the series $\sum_{k=0}^\infty f_k$ converges in 
$L_2$ (and hence in $L_1$) to a function $f$ supported on
$A$ that can be represented as in (\ref{eq4.2}) with $\sup_\alpha
|a_\alpha| \le \delta$. We shall prove that $f$ is a sign. 

Obviously $\int_A f\,d\mu =0$.
Consider $B = \{t \in A \dopu |f(t)| \neq 1 \}$.
By our construction we have for each $n\in\N$
$$
B \subset
\{t \in A \dopu f_n(t) \neq 0 \}
= \Bigl\{t \in A \dopu |f_n(t)| = \frac1m \Bigr\} ,
$$
so $ \mu(B) \le m \|f_n\| $,  and since $\|f_n\| \to 0$, we
conclude that $\mu(B) = 0$.
Therefore $f$ is a sign.
\end{proof}

The previous lemma can also be proved by means of abstract martingale
theory. For simplicity of notation let us work with the classical Haar
system $h_1, h_2 , \dots$ on $[0,1]$. Let $\xi_n= \sum_{k=1}^n h_k$
and $T= \inf\{n\dopu |\xi_n|\le m\}$. Then $(\xi_n)$ is a martingale,
$T$ is a stopping time and $(\xi_n')= (\xi_{n\wedge T})$ is a
uniformly bounded martingale. Hence $(\xi_n')$ converges almost surely
and in $L_1$ to a limit $\xi$ that takes only the values $\pm m$ on
$\{T<\infty\}$. But since $(\xi_n)$ fails to converge pointwise, the
event $\{T=\infty\}$ has probability~$0$. This shows that $\xi=\pm m$
almost surely and $\mathbb{E}\xi=0$. Hence $f=\xi/m$ is the sign we
are looking for.

We are now ready for the main result of this paper. An analogous
theorem for operators on $C(K)$-spaces was proved in \cite{BKSSW}. 

\begin{thm}\label{thm4.3}
Let $T_n \dopu L_1 \to X$ be  HPP-narrow operators, and suppose that
$\sum_{n=1}^\infty T_n$ converges pointwise unconditionally to 
some operator~$T$. Then $T$ is HPP-narrow.
\end{thm}

\begin{proof}
Let $A \in \Sip$, and let $\tilde{\Sig}$ be a nonatomic
sub-$\sigma$-algebra of $\rest{\Sig}{A}$. We have to show
that for every $\eps > 0$ there is a sign 
$f \in L_1(A,\tilde{\Sig}, \mu)$ supported on $A$ with 
$\|Tf\| \le \eps$.

Applying Lemma~\ref{lem3.1} to the restrictions of  $T_n$ and
$T$ to $L_1(A,\tilde{\Sig}, \mu)$ we get a Haar-like 
system $\{h_\alpha \}$ forming a basis for  some 
$L_1^0(A,\Sig_1, \mu)$ and we obtain operators 
$U,V \dopu L_1^0(A,\Sig_1, \mu) \to Y$, 
$W \dopu Y \to X$ such that $\|W\| \le 1$, $T=W(U+V)$ 
on $L_1^0(A,\Sig_1, \mu)$, $\|V\| \le {\eps}/{2}$
and $U$ maps $\{h_\alpha \}$ to a 1-unconditional 
basic sequence. 
Let $C$ be the constant from (\ref{eq4.1}).
 Taking a sign 
$$
f = \sum _{k=0}^\infty
\sum _{\alpha \in \AAA_k}a_\alpha h_\alpha
$$
supported on $A$ with $\sup_\alpha|a_\alpha| \le \eps/(2C)$
(Lemma~\ref{lem4.2}) we obtain from (\ref{eq4.1}) that
$\|Uf\| \le \eps/2$.
Therefore $\|Tf\| \le \|Uf\|+ \|Vf\|\le  \eps$.
\end{proof}

\begin{cor}\label{cor4.3a}
For  any Banach space $X$, no embedding operator is contained in
$\unc (\dots (\unc (\KK (L_1,X))))$.
\end{cor}

\begin{proof}
Compact operators are HPP-narrow.
\end{proof}

The next corollary is due to Rosenthal \cite{Ros-stop}.

\begin{cor}\label{cor4.4}
Every operator $T$ from $L_1$ into a Banach space $X$ with an unconditional
basis is HPP-narrow; in particular it is PP-narrow. Consequently,
$L_1$ does not even sign-embed
into a space with an unconditional basis.
\end{cor}

\begin{proof}
If $P_n$, $n=1,2,\dotsc$, are the partial sum projections associated to
an unconditional basis of $X$, then $T=\sum_{n=1}^\infty (P_n -
P_{n-1}) T $ is a pointwise unconditionally convergent series of
rank-$1$ operators. 
\end{proof}


\section{Questions}

(1) Can one describe $\unc(\cal{K}(L_1, X))$ for general $X$? 
    What about $X=L_1$?

(2) Describe the smallest class of operators $\cal{M} \subset
    \cal{L}(L_1,X)$ that contains the compact operators and is stable
    under pointwise unconditional sums. In particular, is
    $\unc(\cal{K}(L_1,L_1)) = \unc(\unc(\cal{K}(L_1,L_1)))$? 
Note that $X$ does not embed into a space with an unconditional basis
    if  $\cal{M} \neq \cal{L}(L_1,X) $.

(3) Can one develop a similar theory for operators on the James space
    or other spaces that do not embed into spaces with unconditional
    bases?

(4) Is there a space $X$ with the \DP\ such that $\Id \in
    \unc(\dots \allowbreak (\unc(\cal{K}(X,X))))$?

(5) Suppose $E$ is a Banach space with the \DP\ on which the set of
    narrow operators from $E$ to $X$ is a linear space. (This is not
    always the case; e.g., it is not so for $E=X= C([0,1],\ell_1)$
    \cite{BKSSW}.) If $T=\sum T_n$ is a pointwise unconditionally
    convergent series of narrow operators from $E$ into $X$, must $T$
    also be narrow? It is known that under these conditions $\|\Id+T\|
    \ge1$ \cite{KadSSW}. The answer is positive for $E=
    C([0,1],\ell_p)$ if $1<p<\infty$ \cite{BKSSW}.



\end{document}